\documentclass[11pt, oneside, reqno]{amsart}

\usepackage{epsfig}
\usepackage{amsthm}
\usepackage{amssymb}
\usepackage{amsmath}
\usepackage{amscd}
\usepackage{xcolor}
\definecolor{revisionblue}{RGB}{0,76,180}
\usepackage{esint}
\usepackage{enumitem}
\usepackage[top=1.1in, bottom=1in, left=1in, right=1in]{geometry}

\usepackage[hidelinks,pdfstartview=FitV]{hyperref}
\usepackage{mathtools}

\def \<{\langle}
\def \>{\rangle}

\newcommand{\bg}{\begin{equation}}
\newcommand{\ed}{\end{equation}}
\newcommand{\bga}{\begin{eqnarray}}
\newcommand{\eda}{\end{eqnarray}}

\newtheorem {Theorem}  {Theorem}

\numberwithin{Theorem}{section}

\newtheorem {Lemma}[Theorem]  {Lemma}
\newtheorem {Proposition}[Theorem]{Proposition}
\theoremstyle{definition}
\newtheorem{Definition}[Theorem]{Definition}
\theoremstyle{remark}

\expandafter\chardef\csname pre amssym.def
at\endcsname=\the\catcode`\@ \catcode`\@=11
\def\undefine#1{\let#1\undefined}
\def\newsymbol#1#2#3#4#5{\let\next@\relax
 \ifnum#2=\@ne\let\next@\msafam@\else
 \ifnum#2=\tw@\let\next@\msbfam@\fi\fi
 \mathchardef#1="#3\next@#4#5}
\def\mathhexbox@#1#2#3{\relax
 \ifmmode\mathpalette{}{\m@th\mathchar"#1#2#3}%
 \else\leavevmode\hbox{$\m@th\mathchar"#1#2#3$}\fi}
\def\hexnumber@#1{\ifcase#1 0\or 1\or 2\or 3\or 4\or 5\or 6\or 7\or 8\or
 9\or A\or B\or C\or D\or E\or F\fi}

\font\teneufm=eufm10 \font\seveneufm=eufm7 \font\fiveeufm=eufm5
\newfam\eufmfam
\textfont\eufmfam=\teneufm \scriptfont\eufmfam=\seveneufm
\scriptscriptfont\eufmfam=\fiveeufm

\catcode`\@=\csname pre amssym.def at\endcsname

\theoremstyle{remark}
\newtheorem{remark}{Remark}[section]

\numberwithin{equation}{section}
\numberwithin{figure}{section}

\newcommand{\divv}{{\text {div}}\,}

\newcommand{\e}{\epsilon}
\renewcommand{\d}{\delta}
\newcommand{\om}{\omega}

\renewcommand{\k}{\kappa}

\newcommand{\Om}{\Omega}
\renewcommand{\a}{\alpha}

\renewcommand{\b}{\beta}

\newcommand{\gm}{\gamma}
\newcommand{\s}{\sigma}

\newcommand{\Jalpha}{\mathcal{J}_{\alpha}}

\newcommand{\R}{\mathbf{R}}

\newcommand{\les}{\lesssim}

\newcommand{\Ff}{{\mathcal F}}

\newcommand{\Gg}{{\mathcal G}}

\def  \R   {{\mathbb R}}
\def  \Z   {{\mathbb Z}}

\def  \P   {{\mathbb P}}
\def  \T   {{\mathbb T}}

\def  \12  {{\frac{1}{2}}}
\def  \p   {\partial}

\def  \curl  {\nabla \times}

\def\build#1_#2^#3{\mathrel{\mathop{\kern 0pt#1}\limits_{#2}^{#3}}}

\begin{document}

\title{Well-posedness of the relaxed Electron MHD equations with random
diffusion}

\author [Ruimeng Hu]{Ruimeng Hu}
\address{Department of Mathematics and Department of Statistics and Applied Probability, University of California, Santa Barbara,\\ Santa Barbara, CA 93106, USA}
\email{rhu@ucsb.edu} 

\author [Qirui Peng]{Qirui Peng}
\address{Department of Mathematics, University of California, Santa Barbara,\\ Santa Barbara, CA 93106, USA}
\email{qpeng9@ucsb.edu} 

\author [Xu Yang]{Xu Yang}
\address{Department of Mathematics, University of California, Santa Barbara,\\ Santa Barbara, CA 93106, USA}
\email{xuyang@math.ucsb.edu} 

\thanks{}





\begin{abstract}
We investigate a three-dimensional active-vector relaxation of the electron magnetohydrodynamics (EMHD) equations without resistivity. The relaxation replaces the electron current by a generalized current defined through a fractional power of the Laplacian. We drive the system by multiplicative pseudo-differential noise and apply an exponential transformation that removes the stochastic differential term and produces an effective fractional damping. For deterministic divergence-free, mean-zero initial data, we prove local well-posedness up to a stopping time in suitable Gevrey spaces. Under a smallness condition, we further establish global well-posedness with high explicit probability.


\bigskip


\end{abstract}

\maketitle

\section{Introduction}\label{Sec_Intro}

The phenomenon of rapid magnetic reconnection, an essential process in plasma physics responsible for converting magnetic energy into kinetic and thermal energy, is often observed in solar flares, magnetospheric substorms, and laboratory fusion devices. This process involves topological changes in magnetic field lines and the formation of thin current sheets, leading to a rapid release of energy. A widely used mathematical model that captures these dynamics is the Hall-Magnetohydrodynamics (Hall-MHD) system, which extends the classical MHD framework by incorporating the Hall effect to account for the decoupling between electron and ion motions at small spatial scales. The governing equations for the Hall-MHD model take the form:
\begin{subequations}\label{eq:Hall-MHD}
\begin{align}
\partial_t B -\curl (u \times B)+ \nabla \times \big ( (\nabla \times B) \times B \big) &= \eta \Delta B, \label{Hall-MHD 1}\\
\p_t u + u \cdot \nabla u + \nabla p &= (\curl B) \times B + \nu \Delta u, \label{Hall-MHD 2} \\
\nabla \cdot B = 0, \quad \nabla \cdot u &= 0, \label{Hall-MHD 3}
\end{align}
\end{subequations}
where the unknowns \( u, B \), and \( p \) represent the fluid velocity field, magnetic field, and scalar pressure, respectively. The positive constants \( \nu \) and \( \eta \) denote the fluid viscosity and magnetic resistivity; see \cite{B96} for additional physical background.

A closely related model arises by setting \( u \equiv 0 \) in \eqref{eq:Hall-MHD}, which describes the motion of slow ions. This yields the Electron Magnetohydrodynamics (EMHD) equations:
\begin{subequations}\label{eq:EMHD}
\begin{align}
\partial_t B + \nabla \times \big ( (\nabla \times B) \times B \big) &= \eta \Delta B, \label{eq:EMHD_1}\\
\divv B &= 0 ,\label{eq:EMHD_2}\\
B(0) &= B_0.
\end{align}
\end{subequations}

We consider the unknown magnetic field \( B := B(t,x) \) on the domain \( (0,\infty) \times \T^3 \). The strong nonlinearity arising from the term \( \curl \big( (\curl B)\times B\big) \) renders both systems \eqref{eq:Hall-MHD} and \eqref{eq:EMHD} quasilinear. This quasilinear nature introduces significant analytical challenges, making them more difficult to handle than semilinear models like the Navier-Stokes equations. Using the notation \( J = \curl B \) and standard vector identities, equation \eqref{eq:EMHD_1} can be rewritten as:
\begin{equation}\label{eq:EMHD_1_re}
\p_t B + ( B \cdot \nabla)J - (J \cdot \nabla)B = \eta \Delta B.
\end{equation}
A principal analytical difficulty in the EMHD equation \eqref{eq:EMHD_1_re} is the singular term $(B\cdot \nabla)J$.

Recent works \cite{jeong2024illposednesshallelectronmagnetohydrodynamic,dai2024illposedness2frac12delectronmhd} identify strong ill-posedness mechanisms for non-resistive EMHD under the hypotheses considered there. This motivates the search for structurally consistent regularization mechanisms. Following the stochastic Gevrey framework of \cite{BNSW20,HLL25}, we replace deterministic resistivity by multiplicative pseudo-differential noise.

For the nonlinear relaxation, we use the active-vector framework of Chae and Jeong \cite{CJ23}. Their family has the form
\begin{equation}\label{eq:CJ_active_vector_family}
\partial_t B+\curl\big((\curl\mathcal A[B])\times B\big)=0,
\qquad \divv B=0,
\end{equation}
where $\mathcal A$ is a radial Fourier multiplier. We choose
\begin{equation}\label{def:active_vector_multiplier}
\mathcal A_\alpha:=\Lambda^{\alpha-1},
\qquad 0<\alpha\leq1.
\end{equation}
Define the generalized current
\begin{equation}\label{def:active_vector_current}
\Jalpha B:=\curl\Lambda^{\alpha-1}B.
\end{equation}
On the zero Fourier mode we set $\widehat{\Lambda^{\alpha-1}f}(0)=0$, while for $k\neq0$,
\begin{equation}\label{def:active_vector_multiplier_fourier}
\widehat{\Lambda^{\alpha-1}f}(k)=|k|^{\alpha-1}\widehat f(k),
\qquad
\widehat{\Jalpha f}(k)=i k\times |k|^{\alpha-1}\widehat f(k).
\end{equation}
The scalar symbol $|k|^{\alpha-1}$ is real and radial.\\

We consider the stochastic active-vector relaxation
\begin{subequations}\label{eq:EMHD_random_1}
\begin{align}
dB + P_\alpha(B)\,dt &= \mu\Lambda^s B\,dW_t, \label{eq:EMHD_random_1_1}\\
\int_{\T^3}B(x,t)\,dx&=0, \qquad \divv B=0 \label{eq:EMHD_random_1_2}\\
B(0)&=B_0. \label{eq:EMHD_random_1_3}
\end{align}
\end{subequations}
Here $\mu>0$ is the noise amplitude, $0<s\leq1$ and
\begin{equation}\label{def:active_vector_nonlinearity}
P_\alpha(B)
:=\curl\big((\Jalpha B)\times B\big)=(B\cdot\nabla)(\Jalpha B)-(\Jalpha B\cdot\nabla)B.
\end{equation}
The second identity follows because both $B$ and $\Jalpha B$ are divergence-free. The curl form implies
\begin{equation}\label{eq:active_vector_phase_invariance}
\divv P_\alpha(B)=0,
\qquad
\widehat{P_\alpha(B)}(0)=0.
\end{equation}
 When $\alpha=1$, $\Jalpha B=\curl B$ and \eqref{def:active_vector_nonlinearity} is exactly the classical EMHD nonlinearity. For $0<\alpha<1$, $\Jalpha$ has order $\alpha$, so $P_\alpha$ has total differential order $1+\alpha$ and provides a less singular interpolation toward EMHD within the active-vector family of \cite{CJ23}.

We will use below the following Fourier multipliers
\begin{align}
\Ff(\Lambda^s u)(k)&=|k|^s\widehat u(k),\qquad k\in\Z^3, \label{def:Lambda^s}\\
\Ff(e^{\rho\Lambda^s}u)(k)&=e^{\rho|k|^s}\widehat u(k),\qquad \rho>0. \label{def:e^phiLambda^s}
\end{align}
Let $\Gamma=e^{-\mu W_t\Lambda^s}$, $U=\Gamma B$, and
\begin{equation}\label{def:transformed_generalized_current}
V:=\Jalpha U=\Gamma\Jalpha B,
\end{equation}
since $\Jalpha$ commutes with every radial Fourier multiplier. A modewise It\^o calculation gives
\[
d\Gamma=-\mu\Lambda^s\Gamma\,dW_t+\frac12\mu^2\Lambda^{2s}\Gamma\,dt,
\]
and transforms \eqref{eq:EMHD_random_1} into
\begin{subequations}\label{eq:EMHD_Ito_1}
\begin{align}
&\partial_tU+Q(U)=-\frac12\mu^2\Lambda^{2s}U, \label{eq:EMHD_Ito_1_1}\\
&\divv U=0,\quad \int_{\T^3}U(x,t)\,dx=0, \label{eq:EMHD_Ito_1_2}\\
&U(0)=U_0, \label{eq:EMHD_Ito_1_3}
\end{align}
\end{subequations}
where
\begin{align}
Q(U)&:=\Gamma P_\alpha(\Gamma^{-1}U)=Q_1(U)-Q_2(U),
\label{def:transformed_Q}\\
Q_1(U)&:=\Gamma\Big[
(\Gamma^{-1}U\cdot\nabla)(\Gamma^{-1}\Jalpha U)\Big],
\label{def:transformed_Q1}\\
Q_2(U)&:=\Gamma\Big[
(\Gamma^{-1}\Jalpha U\cdot\nabla)(\Gamma^{-1}U)\Big].
\label{def:transformed_Q2}
\end{align}
Equivalently, since $\Gamma$, $\Jalpha$, and $\curl$ commute,
\begin{equation}\label{eq:transformed_Q_curl_form}
Q(U)=\curl\,\Gamma\Big((\Gamma^{-1}\Jalpha U)\times
(\Gamma^{-1}U)\Big).
\end{equation}
For the difference estimate in the fixed-point argument, we introduce the bilinear forms
\begin{subequations}\label{def:transformed_bilinear_Q}
\begin{align}
\mathcal Q_1(U_1,U_2)
&:=\Gamma\Big[(\Gamma^{-1}U_1\cdot\nabla)
(\Gamma^{-1}\Jalpha U_2)\Big],\label{def:transformed_bilinear_Q1}\\
\mathcal Q_2(U_1,U_2)
&:=\Gamma\Big[(\Gamma^{-1}\Jalpha U_1\cdot\nabla)
(\Gamma^{-1}U_2)\Big],\label{def:transformed_bilinear_Q2}\\
\mathcal Q(U_1,U_2)&:=\mathcal Q_1(U_1,U_2)-\mathcal Q_2(U_1,U_2),
\qquad Q(U)=\mathcal Q(U,U).
\label{def:transformed_bilinear_Q_total}
\end{align}
\end{subequations}
The pathwise equation \eqref{eq:EMHD_Ito_1} and the stochastic It\^o equation \eqref{eq:EMHD_random_1} are related rigorously in Section~\ref{sec:probability_space}.

The rest of the paper is organized as follows. In Section~\ref{sec:main_theorem}, we present our main results on local and global well-posedness, along with a brief overview of related work on EMHD and stochastic regularization. Section~\ref{sec:prelim} introduces the necessary notations, functional spaces, and probabilistic setup. In Section~\ref{sec:1D_local_theory}, we prove local well-posedness for the stochastically perturbed EMHD system using stochastic Gevrey estimates. Section~\ref{sec:global} is devoted to the global theory under stronger random diffusion and small initial data. We conclude in Section~\ref{sec:conclusion} with a brief discussion and outlook.

\section{Main theorem and relevant previous works}\label{sec:main_theorem}

Our main results establish local and global well-posedness for the active-vector relaxed system \eqref{eq:EMHD_random_1}. The local lifespan is the minimum of a deterministic contraction time and a Brownian stopping time. The global result holds on an explicit barrier event and requires the Gevrey norm of initial data to be small relative to the effective damping $\mu^2/2$.

\subsection{Main results}

Throughout the statements below, $(\Omega,\mathcal F,(\mathcal F_t)_{t\geq0},\P)$ is a filtered probability space satisfying the usual conditions and carrying a one-dimensional $(\mathcal F_t)$-Brownian motion $W$.

\begin{Theorem}\label{Thm:main_theorem_1}
Let
\begin{equation}\label{eq:parameter_range_main}
\frac34<\a<1,\qquad \frac{1+\a}{2}<s\leq1,
\qquad \frac{1+\a}{s}<\s<2,
\end{equation}
and let $\mu,\a_0,\gm>0$ and $0<\b<\mu^2/2$. Suppose that $B_0\in H\cap\dot{\Gg}^{\s,s}_{\a_0+\gm}$ is deterministic. Define
\begin{equation}\label{eq:local_stopping_main}
T_\om:=\inf\{t\geq0:\mu W_t>\a_0+\b t\}.
\end{equation}
Then there exists a deterministic time
\[
T_*=T_*\!\left(\a,s,\s,\mu,\b,
\|B_0\|_{\dot{\Gg}^{\s,s}_{\a_0+\gm}}\right)>0
\]
such that \eqref{eq:EMHD_random_1} has a unique adapted Gevrey It\^o solution on the stochastic interval
\[
[\![0,\tau_{\rm loc}]\!],\qquad \tau_{\rm loc}:=T_*\wedge T_\om,
\]
and, almost surely,
\[
B\in C\!\left([0,\tau_{\rm loc}];H\cap\dot{\Gg}^{\s,s}_{\gm}\right).
\]
\end{Theorem}

\begin{remark}
The Gevrey weight is used to control the random conjugation
\[
Q(U)=\Gamma P_\alpha(\Gamma^{-1}U).
\]
In view of the similarity between the resulting term and the dissipation on the right-hand side of \eqref{eq:EMHD_Ito_1_1}, it is tempting to take $s$ larger. However, as indicated by the restriction \eqref{eq:parameter_range_main}, $s$ also measures the strength of the noise, and beyond this range the classical PDE estimates can no longer be closed.
\end{remark}

\begin{Theorem}\label{Thm:main_theorem_2}
Assume \eqref{eq:parameter_range_main}. Let $\a_0,\d,\mu>0$, let $0<\b<\mu^2/2$, and let $B_0\in H\cap\dot{\Gg}^{\s,s}_{\a_0+\d}$ be deterministic. There exists a constant $c_{\a,\s,s}>0$ such that, if
\begin{equation}\label{eq:global_smallness_main}
c_{\a,\s,s}\|B_0\|_{\dot{\Gg}^{\s,s}_{\a_0+\d}}<\mu^2-2\b,
\end{equation}
then the unique adapted local solution from Theorem~\ref{Thm:main_theorem_1} extends globally for every sample path in the event below and satisfies
\[
B\in C\!\left([0,\infty);H\cap\dot{\Gg}^{\s,s}_{\d}\right).
\]
The event is
\begin{equation}\label{eq:global_event_main}
\Om^{\a_0,\b}_\mu
:=\{\omega:\mu W_t(\omega)\leq\a_0+\b t\text{ for every }t\geq0\}.
\end{equation}
Moreover,
\begin{equation}\label{eq:global_probability_main}
\P(\Om^{\a_0,\b}_\mu)
=1-\exp\!\left(-\frac{2\a_0\b}{\mu^2}\right).
\end{equation}
In particular, a confidence level $1-\e$ is obtained whenever
$2\a_0\b/\mu^2\geq\log(1/\e)$ and \eqref{eq:global_smallness_main} holds.
\end{Theorem}

\begin{remark}\label{rem:probability_tradeoff}
A convenient specialization makes the probability--damping tradeoff explicit. Fix
$\vartheta\in(0,1/2)$ and take $\b=\vartheta\mu^2$. Then
\[
\P(\Om^{\a_0,\b}_\mu)=1-e^{-2\vartheta\a_0},
\qquad
c_{\a,\s,s}\|B_0\|_{\dot{\Gg}^{\s,s}_{\a_0+\d}}
<(1-2\vartheta)\mu^2.
\]
In particular, the choice $\vartheta=1/4$ and
$\a_0=2\log(1/\e)$ gives probability $1-\e$, while the smallness condition becomes
$\mu^2>2c_{\a,\s,s}\|B_0\|_{\dot{\Gg}^{\s,s}_{\a_0+\d}}$.
\end{remark}

\begin{remark}
The endpoint $s=(1+\a)/2$ is excluded because it would leave no admissible interval for $\s$ under the strict condition $\s>(1+\a)/s$. The classical endpoint $\a=1$ is not addressed by the present estimates.
\end{remark}

\subsection{Related work in the literature}
The study of the well-posedness of Hall-MHD and EMHD systems has received considerable attention over the past decade. For local well-posedness of classical or weak solutions to \eqref{eq:EMHD}, we refer the reader to an incomplete list of works \cite{ADFL11,CDL14,CWW15,D21} and the references therein. More recently, Jeong and Oh \cite{JO25} established local well-posedness of the EMHD system without resistivity under large perturbations around a uniform magnetic field.  Chae and Jeong \cite{CJ23} introduced the active-vector family in \eqref{eq:CJ_active_vector_family}. It contains both three-dimensional Euler and EMHD as distinguished multiplier choices.

Global existence results have also been obtained for reduced models. Dai \cite{dai2023global} proved global existence near a steady state for a $2\frac{1}{2}$-dimensional EMHD system, and later established local well-posedness for the same model with partial resistivity \cite{dai2025well}. A complementary structure-preserving numerical approach for the $2\frac{1}{2}$-dimensional EMHD system with fractional dissipation was developed in \cite{GHPY26}. On the stochastic side, Yamazaki \cite{Y17} constructed martingale solutions for stochastic Hall-MHD equations in both three and $2\frac{1}{2}$ dimensions. More recently, Hu, Peng, and Yang \cite{HPY26} established local pathwise well-posedness and constructed maximal pathwise solutions for the three-dimensional stochastic EMHD system with fractional dissipation and Stratonovich transport noise. By contrast, the active-vector relaxation studied here contains no deterministic resistive term and is driven instead by multiplicative pseudo-differential noise.


Despite various well-posedness results, non-resistive EMHD exhibits strong instability mechanisms. Jeong and Oh \cite{jeong2024illposednesshallelectronmagnetohydrodynamic} obtained ill-posedness and nonexistence results under the geometric and regularity hypotheses specified in their work. Dai showed norm inflation in the $2\frac{1}{2}$-dimensional setting \cite{dai2024illposedness2frac12delectronmhd} and finite-time blow-up for a forced EMHD model in $L^2(\R^3)$ \cite{dai2025blowupforcedelectronmhd}. These results motivate the present active-vector relaxation but should not be read as an unconditional ill-posedness statement for every smooth EMHD datum. For weak solutions, pathological behavior has also been observed. In \cite{P26}, non-uniqueness of steady-state weak solutions in $H^s$ for very low regularity ($s \ll 1$) was established. Moreover, Dai \cite{dai2024nonuniquesolutionselectronmhd} constructed non-unique time-dependent weak solutions in the space $L^1_t L^2_x \cap L^\gamma_t W^{1,\infty}_x$ for any $\gamma < \frac{4}{3}$, further illustrating the ill-posedness of the system in low regularity regimes.


The Gevrey class has emerged as a powerful functional setting for studying analytic and sub-analytic regularity in nonlinear PDEs. Originally introduced by Foias and Temam \cite{FT89} in the context of the Navier–Stokes equations, and later used by Levermore and Oliver \cite{LO97} for the Euler equations, the Gevrey framework has been applied to a wide range of fluid models; see also \cite{KV09,PV11}.

In parallel, stochastic regularization has gained traction as a mechanism to restore well-posedness in otherwise ill-posed PDEs. The work of Buckmaster, Nahmod, Staffilani, and Widmayer \cite{BNSW20} showed high-probability global existence and uniqueness for the surface quasi-geostrophic equation with random diffusion under a Gevrey-type condition. Building on this framework, Hu, Lin, and Liu \cite{HLL25} analyzed the inviscid primitive equations under pseudo-differential noise, proving both local and global well-posedness through delicate stochastic Gevrey estimates.

The present paper develops these ideas for the active-vector relaxed EMHD
system \eqref{eq:EMHD_random_1}. The deterministic nonlinear structure is
inherited from the family introduced in \cite{CJ23}, while the new contribution
is a stochastic Gevrey regularization theory for a nonlinearity of higher
differential order than those treated in \cite{BNSW20,HLL25}. In the surface
quasi-geostrophic and inviscid primitive-equation settings, the relevant
nonlinearities are first-order transport terms. By contrast, $P_\a$ is
quasilinear of order $1+\a$: the generalized current
$\Jalpha=\curl\Lambda^{\a-1}$ contributes $\a$ derivatives in addition to
the transport derivative. After exponential conjugation, the two transport
terms carry the asymmetric Fourier weights $|\ell|^{1+\a}$ and
$|j|^\a|\ell|$, as recorded in
\eqref{eq:K1_Fourier_bound}--\eqref{eq:K2_Fourier_bound}. Their different
frequency placements are controlled by the mixed convolution estimate of
Lemma~\ref{lemma:EMHD_nonlinear_L2}.

The parameter restrictions reflect the two parts of the argument. The local
semigroup estimate requires $\s s>1+\a$ and produces the time singularity
$(t-r)^{-\s/2}$, whose integrability requires $\s<2$. The nonlinear energy
estimate requires $2s>1+\a$, while the assumption $\a>3/4$ provides the
lattice-summability margins used in the convolution bounds. Within this
parameter range, we establish stochastic Gevrey regularization, a stopped
local theory, and high-probability global continuation on the Brownian
barrier event.

Recent work in stochastic PDEs has significantly improved our understanding of noise-induced regularization and long-time dynamics in fluid systems. Several contributions by Hu and collaborators fall into this direction. These include well-posedness results for fractionally dissipated primitive equations with transport noise \cite{HLL_FracPE25} and the regularization of inviscid primitive equations via stochastic Gevrey methods \cite{HLL25}. Stochastic electrohydrodynamics has been explored through the Nernst–Planck–Boussinesq and Nernst–Planck–Navier–Stokes systems, with results on pathwise solutions and ergodicity \cite{HAL_NPB24, HAL_NPNS23}. Additional studies address anisotropic noise estimation in stochastic primitive equations \cite{HCL_AnisoPE25} and pathwise well-posedness for hydrostatic Euler and Navier–Stokes equations under the local Rayleigh condition \cite{HL_HydroEuler23}. For the stochastic inviscid primitive equations, both martingale solutions and pathwise uniqueness have been established \cite{HL_StochPE23}. Together, these works underscore the effectiveness of stochastic modeling in stabilizing complex fluid dynamics, providing motivation for its application to the EMHD system studied here.

\section{Preliminaries}\label{sec:prelim}
\subsection{Notations}
We write $X\les Y$ if $X\leq CY$ for a constant $C>0$. A subscript indicates the parameters on which the constant may depend. We use the normalized Fourier convention on
\[
\T^3=(\R/2\pi\Z)^3,
\qquad
u(x)=\sum_{k\in\Z^3}\widehat u(k)e^{ik\cdot x},
\qquad
\widehat u(k)=\frac{1}{(2\pi)^3}\int_{\T^3}u(x)e^{-ik\cdot x}\,dx.
\]
For vector fields, the Fourier transform is taken componentwise.

\subsection{Sobolev spaces, Gevrey classes, and the phase space}
Let $u:\T^3\to\R^3$ be a periodic vector field. With the normalized measure on $\T^3$, we write
\begin{equation}\label{def:Lp_norm}
\|u\|_{L^p(\T^3)}
:=\left(\frac{1}{(2\pi)^3}\int_{\T^3}|u(x)|^p\,dx\right)^{1/p},
\qquad 1\leq p<\infty,
\end{equation}
and use the usual essential supremum when $p=\infty$. The normalized $L^2$ inner product is
\begin{equation}\label{def:L2_inner_product}
\langle u,v\rangle
:=\frac{1}{(2\pi)^3}\int_{\T^3}u(x)\cdot v(x)\,dx.
\end{equation}
Under the Fourier convention in Section~\ref{sec:prelim}, Parseval's identity becomes
\[
\langle u,v\rangle
=\sum_{k\in\Z^3}\widehat u(k)\cdot\overline{\widehat v(k)},
\qquad
\|u\|_{L^2}^2=\sum_{k\in\Z^3}|\widehat u(k)|^2.
\]

For $r\geq0$, define
\begin{align}
\|u\|_{H^r}^2
&:=\sum_{k\in\Z^3}(1+|k|^2)^r|\widehat u(k)|^2,
\label{def:inhomogeneous_Sobolev}\\
\|u\|_{\Dot H^r}^2
&:=\sum_{k\in\Z^3\setminus\{0\}}|k|^{2r}|\widehat u(k)|^2.
\label{def:homogeneous_Sobolev}
\end{align}

The phase space used throughout the paper is
\begin{equation}\label{def:phase_space_H}
H:=\left\{u\in L^2(\T^3;\R^3):
\divv u=0,\ \widehat u(0)=0\right\}.
\end{equation}
For $u\in H$, the homogeneous and inhomogeneous Sobolev norms of the same positive order are equivalent.

Let $0<s\leq1$, $\s>0$, and $\rho\geq0$. The inhomogeneous Gevrey norm is
\begin{equation}\label{def:inhomogeneous_gevrey_norm}
\begin{aligned}
\|u\|_{\Gg^{\s,s}_{\rho}}^2
&:=\|u\|_{L^2}^2
 +\|e^{\rho\Lambda^s}u\|_{\Dot H^{\s s}}^2\\
&=\sum_{k\in\Z^3}
\left(1+e^{2\rho|k|^s}|k|^{2\s s}\right)|\widehat u(k)|^2,
\end{aligned}
\end{equation}
with the convention that the second summand vanishes at $k=0$. We set
\[
\Gg^{\s,s}_{\rho}
:=\{u\in H^{\s s}:\|u\|_{\Gg^{\s,s}_{\rho}}<\infty\}.
\]
For mean-zero fields it is equivalent, and more convenient, to use the homogeneous norm
\begin{equation}\label{def:gevrey_norm}
\|u\|_{\dot{\Gg}^{\s,s}_{\rho}}^2
:=\sum_{k\in\Z^3\setminus\{0\}}
|k|^{2\s s}e^{2\rho|k|^s}|\widehat u(k)|^2
=\|e^{\rho\Lambda^s}u\|_{\Dot H^{\s s}}^2.
\end{equation}
The corresponding homogeneous Gevrey class is
\[
\dot{\Gg}^{\s,s}_{\rho}
:=\{u\in\Dot H^{\s s}:\|u\|_{\dot{\Gg}^{\s,s}_{\rho}}<\infty\}.
\]
In particular, on $H$,
\begin{equation}\label{eq:homogeneous_inhomogeneous_gevrey_equivalence}
\|u\|_{\Gg^{\s,s}_{\rho}}
\sim_{\s,s,\rho}\|u\|_{\dot{\Gg}^{\s,s}_{\rho}}.
\end{equation}

For a continuous radius $\rho:[0,T]\to[0,\infty)$, define
\begin{equation}\label{def:CTGevrey_space}
C_T\dot{\Gg}^{\s,s}_{\rho}
:=\Big\{u\in C([0,T];H):\;
 e^{\rho(\cdot)\Lambda^s}u(\cdot)
 \in C([0,T];\Dot H^{\s s}),\;
\sup_{0\leq t\leq T}
\|u(t)\|_{\dot{\Gg}^{\s,s}_{\rho(t)}}<\infty\Big\},
\end{equation}
with the norm
\begin{equation}\label{def:CTGevrey_norm}
\|u\|_{C_T\dot{\Gg}^{\s,s}_{\rho}}
:=\sup_{0\leq t\leq T}
\|u(t)\|_{\dot{\Gg}^{\s,s}_{\rho(t)}}.
\end{equation}
This formulation makes explicit that continuity is measured after applying the time-dependent Gevrey weight.

Throughout the paper we use the affine radius
\begin{equation}\label{def:gevrey_radius_phi}
\phi(t):=\a_0+\b t,
\qquad \a_0>0,
\qquad 0<\b<\frac12\mu^2.
\end{equation}
Here $\mu$ is the noise amplitude; the effective damping coefficient in the transformed equation is $\mu^2/2$.

The scalar operator $\Lambda^{\alpha-1}$ is a real radial Fourier multiplier, while
$\Jalpha=\curl\Lambda^{\alpha-1}$ is a real-valued, rotationally covariant multiplier of order $\alpha$. It satisfies
\begin{equation}\label{eq:active_vector_multiplier_properties}
[\Jalpha,\Lambda^r]=[\Jalpha,e^{\rho\Lambda^s}]=0,
\qquad
\divv(\Jalpha u)=0,
\qquad
\widehat{\Jalpha u}(0)=0.
\end{equation}
Moreover, for every periodic vector field $F$,
\begin{equation}\label{eq:curl_phase_properties}
\divv(\curl F)=0,
\qquad
\widehat{\curl F}(0)=0.
\end{equation}
Consequently, $P_\alpha(u)$ and $Q(u)$ belong to $H$ whenever the corresponding expressions are defined and $u\in H$.

We will repeatedly use the following elementary estimate for trading Gevrey radius for derivatives.

\begin{Lemma} \label{lemma:radius_loss}
Let $m\geq0$, $\rho>\k>0$, and
$u\in\dot{\Gg}^{\s+m,s}_{\rho}$. Then
\begin{equation}\label{eq:radius_loss}
\|u\|_{\dot{\Gg}^{\s+m+1,s}_{\rho-\k}}
\leq C_s\k^{-1}
\|u\|_{\dot{\Gg}^{\s+m,s}_{\rho}}.
\end{equation}
In particular, continuity in
$\dot{\Gg}^{\s,s}_{\rho(t)}$ implies continuity in
$\dot{\Gg}^{\s+1,s}_{\rho(t)-\k}$ whenever the radius gap $\k$ is fixed and positive.
\end{Lemma}
\begin{proof}
For $y\geq0$, $ye^{-\k y}\leq C\k^{-1}$. Taking $y=|k|^s$ gives
$|k|^se^{-\k|k|^s}\leq C_s\k^{-1}$.
Multiplying this bound by
$|k|^{(\s+m)s}e^{\rho|k|^s}|\widehat u(k)|$
and summing over $k\neq0$ proves \eqref{eq:radius_loss}.
\end{proof}

\subsection{Probability space and solution concept}\label{sec:probability_space}

Let $(\Omega,\mathcal F,(\mathcal F_t)_{t\geq0},\P)$ satisfy the usual conditions and carry a one-dimensional $(\mathcal F_t)$-Brownian motion $W$. Equivalently, one may use the canonical Wiener space
\[
\Omega=\{\omega\in C([0,\infty);\R):\omega(0)=0\}
\]
with its usual augmented filtration. For a stopping time $\tau$, we use the stochastic-interval notation
\begin{equation}\label{def:stochastic_interval}
[\![0,\tau]\!]
:=\{(\omega,t)\in\Omega\times[0,\infty):0\leq t\leq\tau(\omega)\}.
\end{equation}
Define
\begin{equation}\label{def:Stopping_time}
T_\om:=\inf\{t\geq0:\mu W_t>\phi(t)\}.
\end{equation}
This is a stopping time and $T_\om>0$ almost surely. By continuity of Brownian paths,
\begin{equation}\label{bound:noise_coefficient}
\phi(t)-\mu W_t\geq0,
\qquad 0\leq t\leq T_\om.
\end{equation}
It is useful for the local construction to introduce the stopped path
$\widetilde W_t:=W_{t\wedge T_\om}$. Since $\phi$ is increasing and
$\mu W_{T_\om}=\phi(T_\om)$ on $\{T_\om<\infty\}$, one has
\begin{equation}\label{eq:stopped_path_barrier}
\phi(t)-\mu\widetilde W_t\geq0,\qquad t\geq0.
\end{equation}
Thus the fixed-point problem may be constructed on the deterministic interval
$[0,T_*]$ using the adapted stopped coefficient and then restricted to
$[\![0,T_*\wedge T_\om]\!]$.
For the global theory, define
\begin{equation}\label{def:global_sample_space}
\Om^{\a_0,\b}_\mu
:=\{\omega:\mu W_t(\omega)\leq\a_0+\b t\text{ for every }t\geq0\}.
\end{equation}

\begin{Lemma}\label{lemma:drift_Brownian_motion}
For $\a_0,\b,\mu>0$,
\begin{equation}\label{lemma:drift_Brownian_motion_1}
\P\!\left(\sup_{t\geq0}\left(W_t-\frac{\b}{\mu}t\right)>\frac{\a_0}{\mu}\right)
=\exp\!\left(-\frac{2\a_0\b}{\mu^2}\right).
\end{equation}
Consequently,
\begin{equation}\label{eq:probability_global_sample_space}
\P(\Om^{\a_0,\b}_\mu)
=1-\exp\!\left(-\frac{2\a_0\b}{\mu^2}\right).
\end{equation}
\end{Lemma}
\begin{proof}
Apply the standard all-time maximum formula for Brownian motion with negative drift; see \cite[Proposition~6.8.1]{R92}.
\end{proof}

Next, we recall the definition of weak and mild solutions.
\begin{Definition}\label{def:adapted_solution}
Let $\tau$ be a stopping time and $\gm\geq0$. An
$H\cap\dot{\Gg}^{\s,s}_{\gm}$-valued adapted process $B$ is an adapted
Gevrey It\^o solution of \eqref{eq:EMHD_random_1} on
$[\![0,\tau]\!]$ with radius $\gm$ if the following properties hold.
\begin{enumerate}[label=(\roman*),leftmargin=2.2em]
\item For almost every $\omega\in\Omega$,
\begin{equation}\label{eq:pathwise_continuity_solution}
B(\cdot,\omega)\in
C\!\left([0,\tau(\omega)];
H\cap\dot{\Gg}^{\s,s}_{\gm}\right).
\end{equation}
\item Almost surely,
\begin{equation}\label{eq:solution_integrability}
P_\alpha(B)\in L^1_{\mathrm{loc}}([0,\tau];L^2),
\qquad
\Lambda^sB\in L^2_{\mathrm{loc}}([0,\tau];L^2).
\end{equation}
\item For every test field
$\psi\in C^\infty(\T^3;\R^3)\cap H$ and every $t\geq0$,
\begin{align}\label{eq:weak_Ito_solution}
\langle B(t\wedge\tau),\psi\rangle
&=\langle B_0,\psi\rangle
-\int_0^{t\wedge\tau}
  \langle P_\alpha(B(r)),\psi\rangle\,dr\\
&\quad+\mu\int_0^{t\wedge\tau}
  \langle\Lambda^sB(r),\psi\rangle\,dW_r
\notag
\end{align}
holds almost surely. 
\item The transformed process
$U=e^{-\mu W_t\Lambda^s}B$ satisfies, almost surely,
\begin{equation}\label{eq:solution_transformed_class}
U\in C\!\left([0,\tau];H\cap
\dot{\Gg}^{\s,s}_{\phi(t)+\gm}\right).
\end{equation}
\end{enumerate}
The pairings are taken with respect to the normalized inner product
\eqref{def:L2_inner_product}. 
\end{Definition}

A pathwise mild solution of \eqref{eq:EMHD_Ito_1} on $[t_0,T]$ satisfies
\begin{equation}\label{eq:pathwise_mild_solution}
U(t)=e^{-\frac12\mu^2(t-t_0)\Lambda^{2s}} U(t_0)-\int_{t_0}^te^{-\frac12\mu^2(t-r)\Lambda^{2s}} Q(U(r))\,dr.
\end{equation}

\begin{Proposition}\label{prop:transform_equivalence}
Let $\tau$ be a stopping time and $\gm\geq0$, and assume that
$\phi(t)-\mu W_t\geq0$ on $[\![0,\tau]\!]$. Suppose that $U$ is adapted,
satisfies \eqref{eq:pathwise_mild_solution} on every compact subinterval of
$[\![0,\tau]\!]$, and, almost surely,
\[
U\in C\!\left([0,\tau];H\cap
\dot{\Gg}^{\s,s}_{\phi(t)+\gm}\right).
\]
Then $B=e^{\mu W_t\Lambda^s}U$ is an adapted Gevrey It\^o solution with
radius $\gm$ in the sense of Definition~\ref{def:adapted_solution}.
Conversely, an adapted Gevrey It\^o solution in that class gives a
pathwise mild solution $U=e^{-\mu W_t\Lambda^s}B$ of
\eqref{eq:EMHD_Ito_1}. The transformation is invertible and preserves
uniqueness within these corresponding classes.
\end{Proposition}

\begin{proof}
Fix $k\in\Z^3\setminus\{0\}$. The mild formulation implies that
$t\mapsto\widehat U(k,t)$ is absolutely continuous on compact time intervals and satisfies
\begin{equation}\label{eq:modewise_transformed_equation}
 d\widehat U(k)
 =-\widehat Q(k)\,dt-\frac12\mu^2|k|^{2s}\widehat U(k)\,dt.
\end{equation}
Set
$\widehat B(k)=e^{\mu W_t|k|^s}\widehat U(k)$. Since
$\widehat U(k)$ has finite variation, its quadratic covariation with
$e^{\mu W_t|k|^s}$ vanishes. The scalar It\^o product rule and
\eqref{eq:modewise_transformed_equation} therefore give
\begin{align*}
 d\widehat B(k)
 &=e^{\mu W_t|k|^s}\,d\widehat U(k)
 +\mu|k|^s e^{\mu W_t|k|^s}\widehat U(k)\,dW_t
 +\frac12\mu^2|k|^{2s}e^{\mu W_t|k|^s}\widehat U(k)\,dt\\
 &=-e^{\mu W_t|k|^s}\widehat Q(k)\,dt
 +\mu|k|^s\widehat B(k)\,dW_t.
\end{align*}
By the definition $Q(U)=\Gamma P_\alpha(\Gamma^{-1}U)$,
\[
e^{\mu W_t|k|^s}\widehat Q(k)=\widehat{P_\alpha(B)}(k),
\]
and hence
\begin{equation}\label{eq:modewise_original_equation}
 d\widehat B(k)
 =-\widehat{P_\alpha(B)}(k)\,dt
 +\mu|k|^s\widehat B(k)\,dW_t.
\end{equation}
The zero mode remains zero because both the initial datum and the nonlinear term have zero mean.

Let $\psi\in C^\infty(\T^3;\R^3)\cap H$ and let $\Pi_N$ be the Fourier projection onto
$0<|k|\leq N$. Summing \eqref{eq:modewise_original_equation} against
$\Pi_N\psi$ gives \eqref{eq:weak_Ito_solution} with $\Pi_N\psi$ in place of $\psi$.
Since $\Pi_N\psi\to\psi$ in every Sobolev space, the deterministic pairings converge by the local estimate
$P_\alpha(B)\in L^1_{\rm loc}([0,\tau];L^2)$, which follows from the same convolution bounds used in Section~\ref{sec:1D_local_theory}. Moreover,
$\Lambda^sB\in L^2_{\rm loc}([0,\tau];L^2)$ because $\s s>s$ and $B$ has the stated Gevrey regularity. After localization at the stopping times
\[
\tau_m:=\tau\wedge\inf\left\{t\geq0:
\int_0^t\|\Lambda^sB(r)\|_{L^2}^2\,dr>m\right\},
\]
the scalar It\^o isometry yields convergence of the stochastic pairings. Letting first $N\to\infty$ and then $m\to\infty$ proves \eqref{eq:weak_Ito_solution} for every admissible smooth test field.

Conversely, take trigonometric test fields in Definition~\ref{def:adapted_solution} to recover the modewise equation \eqref{eq:modewise_original_equation}. Applying It\^o's formula to
$e^{-\mu W_t|k|^s}\widehat B(k)$ gives \eqref{eq:modewise_transformed_equation}; summation and using Duhamel's formula yield \eqref{eq:pathwise_mild_solution}. Since the exponential multiplier at time $t$ depends only on $W_t$, the transformation preserves adaptedness. Finally, the exponential transform is invertible on every stopped interval under consideration, and therefore uniqueness in either formulation implies uniqueness in the other.
\end{proof}

\section{Proof of Theorem \ref{Thm:main_theorem_1}}\label{sec:1D_local_theory}
We first prove a time-translated pathwise local result for the transformed equation. Its uniformity with respect to the initial time is later used in the global continuation argument.

\begin{Proposition}\label{prop:EMHD_main_prop_1}
Assume \eqref{eq:parameter_range_main} and $0<\b<\mu^2/2$. Let $t_0\geq0$, $\rho_0>0$, and
\[
\rho(t):=\rho_0+\b(t-t_0).
\]
Fix a continuous sample path and suppose that
\begin{equation}\label{eq:local_path_barrier}
\rho(t)-\mu W_t\geq0
\end{equation}
throughout the time interval under consideration. For every $U_*\in H\cap\dot{\Gg}^{\s,s}_{\rho_0}$ there exists
\[
T_*=T_*\!\left(\a,s,\s,\mu,\b,
\|U_*\|_{\dot{\Gg}^{\s,s}_{\rho_0}}\right)>0,
\]
such that \eqref{eq:EMHD_Ito_1}, with $U(t_0)=U_*$, has a unique mild solution in
\[
C\!\left([t_0,t_0+T];H\cap\dot{\Gg}^{\s,s}_{\rho(t)}\right)
\]
for every $0<T\leq T_*$ on which \eqref{eq:local_path_barrier} holds. If
$\|U_*\|_{\dot{\Gg}^{\s,s}_{\rho_0}}\leq M$, the same positive time $T_*(M)$ may be chosen for all $t_0$, all such paths, and all such data.
\end{Proposition}

\begin{proof}
Let
\[
S(t):=e^{-\frac12\mu^2t\Lambda^{2s}},
\qquad
M:=\|U_*\|_{\dot{\Gg}^{\s,s}_{\rho_0}},
\]
and, for $T>0$, set
\[
X_{t_0,T}:=C\!\left([t_0,t_0+T];H\cap
\dot{\Gg}^{\s,s}_{\rho(t)}\right),
\qquad
\|U\|_{X_{t_0,T}}
:=\sup_{t_0\leq t\leq t_0+T}
\|U(t)\|_{\dot{\Gg}^{\s,s}_{\rho(t)}}.
\]
Define the Duhamel map
\begin{equation}\label{def:sol_operator}
\Phi(U)(t):=S(t-t_0)U_*-\int_{t_0}^tS(t-r)Q(U(r))\,dr.
\end{equation}
The semigroup preserves $H$, while $Q(U)$ is a curl and has zero Fourier mode. Thus $\Phi$ maps $H$ into itself.

We first record the two transformed nonlinearities in Fourier variables. Put
$\ell=k-j$. From \eqref{def:transformed_Q1}--\eqref{def:transformed_Q2},
\begin{align}
\widehat{Q_1(U)}(k,r)
&=e^{-\mu W_r|k|^s}
\sum_{j+\ell=k}
 e^{\mu W_r|j|^s}
 \big(\widehat U(j,r)\cdot i\ell\big)
 e^{\mu W_r|\ell|^s}
 \big(i\ell\times |\ell|^{\a-1}\widehat U(\ell,r)\big),
\label{eq:Fourier_Q1_exact}\\
\widehat{Q_2(U)}(k,r)
&=e^{-\mu W_r|k|^s}
\sum_{j+\ell=k}
 e^{\mu W_r|j|^s}
 \big(i j\times |j|^{\a-1}\widehat U(j,r)\big)\cdot i\ell
 e^{\mu W_r|\ell|^s}\widehat U(\ell,r).
\label{eq:Fourier_Q2_exact}
\end{align}
All sums may be restricted to nonzero input modes because $U\in H$. In particular,
\begin{align}
|\widehat{Q_1(U)}(k,r)|
&\les_\a e^{-\mu W_r|k|^s}
\sum_{j+\ell=k}
 e^{\mu W_r|j|^s}|\widehat U(j,r)|
 e^{\mu W_r|\ell|^s}|\ell|^{1+\a}|\widehat U(\ell,r)|,
\label{eq:Fourier_Q1_abs}\\
|\widehat{Q_2(U)}(k,r)|
&\les_\a e^{-\mu W_r|k|^s}
\sum_{j+\ell=k}
 e^{\mu W_r|j|^s}|j|^\a|\widehat U(j,r)|
 e^{\mu W_r|\ell|^s}|\ell|\,|\widehat U(\ell,r)|.
\label{eq:Fourier_Q2_abs}
\end{align}

Let $t_0\leq r<t\leq t_0+T$ and put $h:=t-r$. Since
$\rho(t)-\rho(r)=\b h$, the subadditivity
\begin{equation}\label{bound:Fourier_derivatives}
|k|^s\leq|j|^s+|\ell|^s,
\qquad 0<s\leq1,
\end{equation}
and the path barrier $\rho(r)-\mu W_r\geq0$ imply
\begin{align}
&e^{\rho(t)|k|^s}e^{-\mu W_r|k|^s}
 e^{\mu W_r|j|^s}e^{\mu W_r|\ell|^s}\notag\\
&\qquad
=e^{\b h|k|^s}
 e^{(\rho(r)-\mu W_r)|k|^s}
 e^{\mu W_r(|j|^s+|\ell|^s)}\notag\\
&\qquad
\leq e^{\b h|k|^s}
 e^{\rho(r)|j|^s}e^{\rho(r)|\ell|^s}.
\label{eq:random_weight_split}
\end{align}
Since $|k|\geq1$ on $H$ and $\b<\mu^2/2$, one has 
\begin{align}
|k|^{\s s}e^{\b h|k|^s-\frac12\mu^2h|k|^{2s}}
&\leq |k|^{\s s}e^{-(\frac12\mu^2-\b)h|k|^{2s}}
\les_{\mu,\b,\s,s} h^{-\s/2}.
\label{eq:propagator_est_1}
\end{align}
Indeed, this follows after setting $y=h^{1/2}|k|^s$ and using
$\sup_{y\geq0}y^\s e^{-c y^2}<\infty$.

We now estimate the two Duhamel integrands separately. Define
\begin{align*}
I_1(t,r)&:=\left\|e^{\rho(t)\Lambda^s}
\Lambda^{\s s}S(t-r)Q_1(U(r))\right\|_{L^2},\\
I_2(t,r)&:=\left\|e^{\rho(t)\Lambda^s}
\Lambda^{\s s}S(t-r)Q_2(U(r))\right\|_{L^2}.
\end{align*}
By \eqref{eq:Fourier_Q1_abs}, \eqref{eq:random_weight_split}, and
\eqref{eq:propagator_est_1},
\begin{align}
I_1(t,r)
&\les h^{-\s/2}
\bigg(\sum_{k\neq0}
\Big(\sum_{j+\ell=k}
 e^{\rho(r)|j|^s}|\widehat U(j,r)|
 e^{\rho(r)|\ell|^s}|\ell|^{1+\a}|\widehat U(\ell,r)|
\Big)^2\bigg)^{1/2}\notag\\
&\les h^{-\s/2}
\Big(\sum_{j\neq0}e^{\rho(r)|j|^s}|\widehat U(j,r)|\Big)
\Big(\sum_{\ell\neq0}e^{2\rho(r)|\ell|^s}
 |\ell|^{2(1+\a)}|\widehat U(\ell,r)|^2\Big)^{1/2}\notag\\
&\les h^{-\s/2}
\Big(\sum_{j\neq0}|j|^{-2\s s}\Big)^{1/2}
\|U(r)\|_{\dot{\Gg}^{\s,s}_{\rho(r)}}
\Big(\sum_{\ell\neq0}e^{2\rho(r)|\ell|^s}
 |\ell|^{2\s s}|\widehat U(\ell,r)|^2\Big)^{1/2}\notag\\
&\les_{\a,\s,s}h^{-\s/2}
\|U(r)\|_{\dot{\Gg}^{\s,s}_{\rho(r)}}^2.
\label{eq:local_I1_detail}
\end{align}
Here we used Young's inequality $\ell^1*\ell^2\to\ell^2$,
$2\s s>3$, and $1+\a<\s s$.

For the second term, Young's inequality
$\ell^{4/3}*\ell^{4/3}\to\ell^2$ gives
\begin{align}
I_2(t,r)
&\les h^{-\s/2}
\left\|e^{\rho(r)|j|^s}|j|^\a\widehat U(j,r)\right\|_{\ell^{4/3}}
\left\|e^{\rho(r)|\ell|^s}|\ell|\widehat U(\ell,r)\right\|_{\ell^{4/3}}.
\label{eq:local_I2_young}
\end{align}
H\"older's inequality yields
\begin{align}
\left\|e^{\rho(r)|j|^s}|j|^\a\widehat U(j,r)\right\|_{\ell^{4/3}}
&\leq
\|U(r)\|_{\dot{\Gg}^{\s,s}_{\rho(r)}}
\Big(\sum_{j\neq0}|j|^{-4(\s s-\a)}\Big)^{1/4},
\label{eq:local_Q2_factor1}\\
\left\|e^{\rho(r)|\ell|^s}|\ell|\widehat U(\ell,r)\right\|_{\ell^{4/3}}
&\leq
\|U(r)\|_{\dot{\Gg}^{\s,s}_{\rho(r)}}
\Big(\sum_{\ell\neq0}|\ell|^{-4(\s s-1)}\Big)^{1/4}.
\label{eq:local_Q2_factor2}
\end{align}
The two lattice sums converge because
$\s s-\a>1>3/4$ and $\s s-1>\a>3/4$. Consequently,
\begin{equation}\label{eq:local_I2_detail}
I_2(t,r)\les_{\a,\s,s}h^{-\s/2}
\|U(r)\|_{\dot{\Gg}^{\s,s}_{\rho(r)}}^2.
\end{equation}
Combining \eqref{eq:local_I1_detail} and \eqref{eq:local_I2_detail},
\begin{equation}\label{eq:semigroup_Q_local}
\left\|S(t-r)Q(U(r))\right\|_{\dot{\Gg}^{\s,s}_{\rho(t)}}
\les (t-r)^{-\s/2}
\|U(r)\|_{\dot{\Gg}^{\s,s}_{\rho(r)}}^2.
\end{equation}

For the linear term, $|k|^s\leq|k|^{2s}$ on nonzero modes gives
\begin{equation}\label{eq:local_linear_term}
\|S(t-t_0)U_*\|_{\dot{\Gg}^{\s,s}_{\rho(t)}}
\leq\|U_*\|_{\dot{\Gg}^{\s,s}_{\rho_0}}=M.
\end{equation}

Let $\theta:=1-\s/2>0$. Integrating
\eqref{eq:semigroup_Q_local} and using \eqref{eq:local_linear_term}, we obtain
\begin{equation}\label{est:1D_solution_map_3}
\|\Phi(U)\|_{X_{t_0,T}}
\leq M+C T^\theta\|U\|_{X_{t_0,T}}^2.
\end{equation}
To estimate differences, use the bilinear forms in
\eqref{def:transformed_bilinear_Q}. For example,
\[
Q_1(U)-Q_1(\widetilde U)
=\mathcal Q_1(U-\widetilde U,U)
+\mathcal Q_1(\widetilde U,U-\widetilde U),
\]
and analogously for $Q_2$. Repeating the preceding estimates gives
\begin{equation}\label{est:1D_solution_map_4}
\|\Phi(U)-\Phi(\widetilde U)\|_{X_{t_0,T}}
\leq C T^\theta
\big(\|U\|_{X_{t_0,T}}+\|\widetilde U\|_{X_{t_0,T}}\big)
\|U-\widetilde U\|_{X_{t_0,T}}.
\end{equation}

If $M=0$, the zero solution is unique. Otherwise, consider the closed ball
$\mathbb B_{2M}\subset X_{t_0,T}$. Choose $T_*(M)>0$ so that
\begin{equation}\label{eq:contraction_time_choice}
4CM T_*(M)^\theta\leq\frac12.
\end{equation}
Then \eqref{est:1D_solution_map_3} maps $\mathbb B_{2M}$ into itself and
\eqref{est:1D_solution_map_4} makes $\Phi$ a strict contraction there. This proves existence and uniqueness. All constants depend only on the fixed parameters and $M$, not on $t_0$, the particular path satisfying the barrier, or the value of the initial radius. The same estimates applied to the Picard iterates show that an $\mathcal F_{t_0}$-measurable restart datum produces an adapted solution when the coefficient path is adapted.
\end{proof}

\subsection{Proof of Theorem \ref{Thm:main_theorem_1}}\label{sec:proof_main_thm_1}

Set $\rho(t)=\phi(t)+\gm$ and use the stopped path
$\widetilde W_t=W_{t\wedge T_\om}$. By \eqref{eq:stopped_path_barrier},
$\rho(t)-\mu\widetilde W_t\geq\gm$ on the deterministic interval $[0,T_*]$.
Proposition~\ref{prop:EMHD_main_prop_1} therefore gives a unique solution there.
The Picard iterates are adapted because they are formed from the adapted stopped path, and their uniform limit is adapted. Restricting to
$\tau_{\rm loc}=T_*\wedge T_\om$, where $\widetilde W=W$, gives
\[
U\in C\!\left([0,\tau_{\rm loc}];H\cap
\dot{\Gg}^{\s,s}_{\phi(t)+\gm}\right).
\]
Let $B=e^{\mu W_t\Lambda^s}U$. Then
\begin{equation}\label{eq:B_from_U_local_bound}
\|B(t)\|_{\dot{\Gg}^{\s,s}_{\gm}}
=\|e^{-(\phi(t)-\mu W_t)\Lambda^s}
 e^{(\phi(t)+\gm)\Lambda^s}\Lambda^{\s s}U(t)\|_{L^2}
\leq\|U(t)\|_{\dot{\Gg}^{\s,s}_{\phi(t)+\gm}}.
\end{equation}
The function $t\mapsto\phi(t)-\mu W_t$ is continuous and nonnegative on the stopped interval. Hence
$
e^{-(\phi(t)-\mu W_t)\Lambda^s}
$
is a strongly continuous family of contractions on $\Dot H^{\s s}$. Combining this fact with \eqref{eq:B_from_U_local_bound} and the weighted continuity of $U$ yields the continuity of $B$ in
$H\cap\dot{\Gg}^{\s,s}_{\gm}$. Proposition~\ref{prop:transform_equivalence} shows that $B$ is an adapted It\^o solution and transfers uniqueness from $U$ to $B$. This proves Theorem~\ref{Thm:main_theorem_1}.

\section{Proof of Theorem \ref{Thm:main_theorem_2}}\label{sec:global}
We establish the global result through a nonlinear energy estimate followed by a radius-loss continuation argument.

\subsection{Energy estimates}
Let $\rho(t)=\rho_0+\b t$ be any affine Gevrey radius satisfying
$\rho(t)-\mu W_t\geq0$ on the interval in question.

\begin{Proposition}\label{prop:EMHD_main_prop_2}
Assume \eqref{eq:parameter_range_main}. Let $U$ solve \eqref{eq:EMHD_Ito_1} and suppose
\[
U\in C\!\left([0,T];H\cap\dot{\Gg}^{\s+1,s}_{\rho(t)}\right).
\]
Then, for almost every $t\in[0,T]$,
\begin{equation}\label{prop:EMHD_main_prop_2_1}
\frac{d}{dt}\|U\|_{\dot{\Gg}^{\s,s}_{\rho(t)}}^2
+\left(\mu^2-2\b-c_{\a,\s,s}
\|U\|_{\dot{\Gg}^{\s,s}_{\rho(t)}}\right)
\|U\|_{\dot{\Gg}^{\s+1,s}_{\rho(t)}}^2\leq0.
\end{equation}
\end{Proposition}

The nonlinear term satisfies the following mixed estimate.
\begin{Lemma}\label{lemma:EMHD_nonlinear_L2}
Under the assumptions of Proposition~\ref{prop:EMHD_main_prop_2},
\begin{equation}\label{lemma:EMHD_nonlinear_L2_1}
\left|\left\langle e^{\rho(t)\Lambda^s}Q(U),
 e^{\rho(t)\Lambda^s}\Lambda^{2\s s}U\right\rangle\right|
\les_{\a,\s,s}
\|U\|_{\dot{\Gg}^{\s,s}_{\rho(t)}}
\|U\|_{\dot{\Gg}^{\s+1,s}_{\rho(t)}}^2.
\end{equation}
\end{Lemma}

\begin{proof}
Fix $t$ and abbreviate
\begin{equation}\label{eq:LH_definitions}
L:=\|U\|_{\dot{\Gg}^{\s,s}_{\rho(t)}},
\qquad
H_1:=\|U\|_{\dot{\Gg}^{\s+1,s}_{\rho(t)}},
\qquad
a_k:=e^{\rho(t)|k|^s}|\widehat U(k)|.
\end{equation}
According to \eqref{def:transformed_Q}--\eqref{def:transformed_Q2}, write
\begin{align}
K&:=\left\langle e^{\rho(t)\Lambda^s}Q(U),
 e^{\rho(t)\Lambda^s}\Lambda^{2\s s}U\right\rangle
=K_1-K_2,
\label{eq:K_total_definition}\\
K_1&:=\left\langle e^{\rho(t)\Lambda^s}Q_1(U),
 e^{\rho(t)\Lambda^s}\Lambda^{2\s s}U\right\rangle,
\label{eq:K1_definition}\\
K_2&:=\left\langle e^{\rho(t)\Lambda^s}Q_2(U),
 e^{\rho(t)\Lambda^s}\Lambda^{2\s s}U\right\rangle.
\label{eq:K2_definition}
\end{align}
Thus $|K|\leq|K_1|+|K_2|$.

For completeness, the two contributions can be written in Fourier space as follows. With
$\ell=k-j$ and harmless normalization constants suppressed,
\begin{align}
K_1
&=\sum_{k\neq0}e^{\rho(t)|k|^s}|k|^{2\s s}
\overline{\widehat U(k)}\cdot
\Big[e^{(\rho(t)-\mu W_t)|k|^s}
\sum_{j+\ell=k}e^{\mu W_t|j|^s}
(\widehat U(j)\cdot i\ell)\notag\\
&\hspace{7.2cm}\times
 e^{\mu W_t|\ell|^s}
 (i\ell\times |\ell|^{\a-1}\widehat U(\ell))\Big],
\label{eq:K1_Fourier_exact}\\
K_2
&=\sum_{k\neq0}e^{\rho(t)|k|^s}|k|^{2\s s}
\overline{\widehat U(k)}\cdot
\Big[e^{(\rho(t)-\mu W_t)|k|^s}
\sum_{j+\ell=k}e^{\mu W_t|j|^s}
\big((ij\times|j|^{\a-1}\widehat U(j))\cdot i\ell\big)\notag\\
&\hspace{7.2cm}\times
 e^{\mu W_t|\ell|^s}\widehat U(\ell)\Big].
\label{eq:K2_Fourier_exact}
\end{align}
Since $\rho(t)-\mu W_t\geq0$, the same subadditivity argument as in
\eqref{eq:random_weight_split} gives
\begin{align}
|K_1|
&\les_\a\sum_{k\neq0}|k|^{2\s s}a_k
\sum_{j+\ell=k}a_j|\ell|^{1+\a}a_\ell,
\label{eq:K1_Fourier_bound}\\
|K_2|
&\les_\a\sum_{k\neq0}|k|^{2\s s}a_k
\sum_{j+\ell=k}|j|^\a a_j|\ell|a_\ell.
\label{eq:K2_Fourier_bound}
\end{align}

Set $q:=(\s-1)s>0$. In \eqref{eq:K1_Fourier_bound}, write
$|k|^{2\s s}=|k|^{(\s+1)s}|k|^q$ and use
\begin{equation}\label{eq:q_split}
|k|^q\les_q |j|^q+|\ell|^q,
\qquad k=j+\ell.
\end{equation}
Cauchy--Schwarz in the output frequency then yields
$|K_1|\les H_1(K_{11}+K_{12})$, where,
\begin{equation}
K_{11}
:=\bigg(\sum_{k\neq0}
\Big(\sum_{j+\ell=k}|j|^q a_j|\ell|^{1+\a}a_\ell\Big)^2\bigg)^{1/2},\quad
K_{12}
:=\bigg(\sum_{k\neq0}
\Big(\sum_{j+\ell=k}a_j|\ell|^{q+1+\a}a_\ell\Big)^2\bigg)^{1/2}.
\label{eq:K12_definition}
\end{equation}
Young's inequality $\ell^{4/3}*\ell^{4/3}\to\ell^2$ and H\"older's inequality give
\begin{align}
K_{11}
&\leq\big\||j|^q a_j\big\|_{\ell^{4/3}}
\big\||j|^{1+\a}a_j\big\|_{\ell^{4/3}}\notag\\
&\leq
\Big(\sum_{j\neq0}|j|^{2\s s}a_j^2\Big)^{1/2}
\Big(\sum_{j\neq0}|j|^{-4s}\Big)^{1/4}\times
\Big(\sum_{j\neq0}|j|^{2(\s+1)s}a_j^2\Big)^{1/2}
\Big(\sum_{j\neq0}|j|^{-4((\s+1)s-(1+\a))}\Big)^{1/4}\notag\\
&\les L H_1.
\label{eq:est_K_11}
\end{align}
For $K_{12}$, Young's inequality $\ell^1*\ell^2\to\ell^2$ yields
\begin{align}
K_{12}
&\leq\|a\|_{\ell^1}
\big\||j|^{q+1+\a}a_j\big\|_{\ell^2}\notag\\
&\leq
\Big(\sum_{j\neq0}|j|^{-2\s s}\Big)^{1/2}L
\big\||j|^{q+1+\a}a_j\big\|_{\ell^2}\les L H_1,
\label{eq:est_K_12}
\end{align}
where the last step uses
$q+1+\a=(\s-1)s+1+\a\leq(\s+1)s$, equivalently
$1+\a\leq2s$. Hence
\begin{equation}\label{eq:K1_final_bound}
|K_1|\les L H_1^2.
\end{equation}

We turn to \eqref{eq:K2_Fourier_bound}. After the same extraction of the output factor and the split \eqref{eq:q_split},
\begin{equation}\label{eq:K2_split}
|K_2|\les H_1(K_{21}+K_{22}),
\end{equation}
where
\begin{equation}
K_{21}
:=\bigg(\sum_{k\neq0}
\Big(\sum_{j+\ell=k}|j|^{q+\a}a_j|\ell|a_\ell\Big)^2\bigg)^{1/2},\quad
K_{22}
:=\bigg(\sum_{k\neq0}
\Big(\sum_{j+\ell=k}|j|^\a a_j|\ell|^{q+1}a_\ell\Big)^2\bigg)^{1/2}.
\end{equation}
Again using $\ell^{4/3}*\ell^{4/3}\to\ell^2$, we obtain
\begin{align}
K_{21}
&\leq\big\||j|^{q+\a}a_j\big\|_{\ell^{4/3}}
\big\||j|a_j\big\|_{\ell^{4/3}}\notag\\
&\leq
H_1\Big(\sum_{j\neq0}|j|^{-4(2s-\a)}\Big)^{1/4}
L\Big(\sum_{j\neq0}|j|^{-4(\s s-1)}\Big)^{1/4}\les L H_1,
\label{eq:est_K_21}\\
K_{22}
&\leq\big\||j|^\a a_j\big\|_{\ell^{4/3}}
\big\||j|^{q+1}a_j\big\|_{\ell^{4/3}}\notag\\
&\leq
L\Big(\sum_{j\neq0}|j|^{-4(\s s-\a)}\Big)^{1/4}
H_1\Big(\sum_{j\neq0}|j|^{-4(2s-1)}\Big)^{1/4}\les L H_1.
\label{eq:est_K_22}
\end{align}
Indeed, the assumptions imply
\begin{align*}
&s>\frac{1+\a}{2}>\frac34,
\qquad \s s>1+\a>\frac32,\\
&(\s+1)s-(1+\a)>s>\frac34,
\qquad 2s-\a>1>\frac34,\\
&\s s-1>\a>\frac34,
\qquad \s s-\a>1>\frac34,
\qquad 2s-1>\a>\frac34.
\end{align*}
Thus every lattice sum in
\eqref{eq:est_K_11}--\eqref{eq:est_K_22} converges. We conclude that
\begin{equation}\label{eq:K2_final_bound}
|K_2|\les L H_1^2.
\end{equation}
Combining \eqref{eq:K1_final_bound} and \eqref{eq:K2_final_bound} proves
\eqref{lemma:EMHD_nonlinear_L2_1}.
\end{proof}

\begin{proof}[Proof of Proposition~\ref{prop:EMHD_main_prop_2}]
By a standard density argument, it is enough to perform the following computation for smooth solutions; the general case follows by spatial regularization, the assumed
$C([0,T];\dot{\Gg}^{\s+1,s}_{\rho(t)})$ regularity, and
Lemma~\ref{lemma:EMHD_nonlinear_L2}. We therefore test
\eqref{eq:EMHD_Ito_1_1} directly against
$e^{2\rho(t)\Lambda^s}\Lambda^{2\s s}U$.

Using that the Fourier multipliers involved are self-adjoint and commute, we obtain
\begin{align}
&\left\langle e^{\rho(t)\Lambda^s}\partial_tU,
 e^{\rho(t)\Lambda^s}\Lambda^{2\s s}U\right\rangle
+\left\langle e^{\rho(t)\Lambda^s}Q(U),
 e^{\rho(t)\Lambda^s}\Lambda^{2\s s}U\right\rangle\notag\\
&\hspace{4cm}=-\frac12\mu^2
\|U\|_{\dot{\Gg}^{\s+1,s}_{\rho(t)}}^2.
\label{eq:tested_energy_identity}
\end{align}
Since $\dot\rho(t)=\b$, the chain rule gives
\begin{equation}\label{eq:weighted_time_derivative_identity}
\left\langle e^{\rho(t)\Lambda^s}\partial_tU,
 e^{\rho(t)\Lambda^s}\Lambda^{2\s s}U\right\rangle
=\frac12\frac{d}{dt}
\|U\|_{\dot{\Gg}^{\s,s}_{\rho(t)}}^2
-\b\|U\|_{\dot{\Gg}^{\s+1/2,s}_{\rho(t)}}^2.
\end{equation}
Combining \eqref{eq:tested_energy_identity} and
\eqref{eq:weighted_time_derivative_identity}, we find for almost every $t\in[0,T]$ that
\begin{align}
\frac12\frac{d}{dt}
\|U\|_{\dot{\Gg}^{\s,s}_{\rho(t)}}^2
-\b\|U\|_{\dot{\Gg}^{\s+1/2,s}_{\rho(t)}}^2
+\left\langle e^{\rho(t)\Lambda^s}Q(U),
 e^{\rho(t)\Lambda^s}\Lambda^{2\s s}U\right\rangle
=-\frac12\mu^2
\|U\|_{\dot{\Gg}^{\s+1,s}_{\rho(t)}}^2.
\label{eq:full_energy_identity}
\end{align}

Since the zero mode vanishes and $|k|\geq1$,
\begin{equation}\label{eq:half_to_one_derivative}
\|U\|_{\dot{\Gg}^{\s+1/2,s}_{\rho(t)}}
\leq\|U\|_{\dot{\Gg}^{\s+1,s}_{\rho(t)}}.
\end{equation}
Applying Lemma~\ref{lemma:EMHD_nonlinear_L2} to
\eqref{eq:full_energy_identity}, using \eqref{eq:half_to_one_derivative}, and multiplying by two gives
\[
\frac{d}{dt}\|U\|_{\dot{\Gg}^{\s,s}_{\rho(t)}}^2
+\left(\mu^2-2\b-c_{\a,\s,s}
\|U\|_{\dot{\Gg}^{\s,s}_{\rho(t)}}\right)
\|U\|_{\dot{\Gg}^{\s+1,s}_{\rho(t)}}^2\leq0,
\]
after enlarging the constant $c_{\a,\s,s}$. This is
\eqref{prop:EMHD_main_prop_2_1}.
\end{proof}

\subsection{Extension of the local solution}

The next proposition combines local existence with radius loss and the energy estimate. Local existence is first obtained with the same $\s$ index at the full radius; a positive radius gap then yields one additional $s$-derivative, allowing the energy estimate to be applied.

\begin{Proposition}\label{prop:EMHD_main_prop_3}
Assume \eqref{eq:parameter_range_main}. Fix $\d>0$, $\a_0>0$, and $0<\b<\mu^2/2$. Suppose
$U_0\in H\cap\dot{\Gg}^{\s,s}_{\a_0+\d}$ and
\begin{equation}\label{prop:EMHD_main_prop_3_1}
c_{\a,\s,s}\|U_0\|_{\dot{\Gg}^{\s,s}_{\a_0+\d}}<\mu^2-2\b.
\end{equation}
Then, for every $\omega\in\Om^{\a_0,\b}_\mu$, system \eqref{eq:EMHD_Ito_1} has a unique global pathwise solution
\[
U\in C\!\left([0,\infty);H\cap
\dot{\Gg}^{\s,s}_{\phi(t)+\d}\right).
\]
Moreover,
\begin{equation}\label{eq:global_U_bound_full_radius}
\|U(t)\|_{\dot{\Gg}^{\s,s}_{\phi(t)+\d}}
\leq\|U_0\|_{\dot{\Gg}^{\s,s}_{\a_0+\d}},
\qquad t\geq0,
\end{equation}
and the left-hand side is nonincreasing in time.
\end{Proposition}

\begin{proof}
Set
\begin{equation}\label{eq:M_global_definition}
M:=\|U_0\|_{\dot{\Gg}^{\s,s}_{\a_0+\d}},
\qquad
\eta:=\mu^2-2\b-c_{\a,\s,s}M>0.
\end{equation}
Fix $0\leq\d'<\d$ and define
\begin{equation}\label{eq:delta_n_sequence}
\d_n:=\d'+2^{-n}(\d-\d'),
\qquad
\rho_n(t):=\phi(t)+\d_n,
\qquad n\geq0.
\end{equation}
Then $\d_0=\d$, $\d_n\downarrow\d'$, and the difference
\begin{equation}\label{eq:kappa_n_definition}
\k_n:=\rho_n(t)-\rho_{n+1}(t)
=2^{-(n+1)}(\d-\d')>0
\end{equation}
is independent of time. Let $T_0:=T_*(M)$ be the common local existence time from
Proposition~\ref{prop:EMHD_main_prop_1}. On the event
$\Om^{\a_0,\b}_\mu$,
\begin{equation}\label{eq:global_barrier_each_radius}
\rho_n(t)-\mu W_t
=\phi(t)-\mu W_t+\d_n\geq\d_n\geq0,
\end{equation}
for all $t\geq0$ and all $n$.

We prove by induction that, for each integer $n\geq1$,
\begin{subequations}\label{eq:global_induction_statements}
\begin{align}
&U\in C\!\left([0,nT_0];H\cap
\dot{\Gg}^{\s,s}_{\rho_n(t)}\right),
\label{eq:global_induction_existence}\\
&t\longmapsto
\|U(t)\|_{\dot{\Gg}^{\s,s}_{\rho_n(t)}}
\text{ is nonincreasing on }[0,nT_0],
\label{eq:global_induction_monotonicity}\\
&\|U(t)\|_{\dot{\Gg}^{\s,s}_{\rho_n(t)}}\leq M,
\qquad 0\leq t\leq nT_0.
\label{eq:global_induction_bound}
\end{align}
\end{subequations}

For the base step, Proposition~\ref{prop:EMHD_main_prop_1} applied at the full radius
$\rho_0(t)=\phi(t)+\d$ gives
\begin{equation}\label{eq:base_local_full_radius}
U\in C\!\left([0,T_0];H\cap
\dot{\Gg}^{\s,s}_{\rho_0(t)}\right).
\end{equation}
The fixed radius gap $\k_0$ and Lemma~\ref{lemma:radius_loss} imply
\begin{equation}\label{eq:base_high_regularity}
U\in C\!\left([0,T_0];H\cap
\dot{\Gg}^{\s+1,s}_{\rho_1(t)}\right).
\end{equation}
To see the continuity explicitly, observe that
\[
e^{\rho_1(t)\Lambda^s}\Lambda^{(\s+1)s}U(t)
=\Lambda^s e^{-\k_0\Lambda^s}
\big(e^{\rho_0(t)\Lambda^s}\Lambda^{\s s}U(t)\big),
\]
and $\Lambda^s e^{-\k_0\Lambda^s}$ is bounded on $L^2$.

We may now apply Proposition~\ref{prop:EMHD_main_prop_2} at the radius $\rho_1(t)$. Since
$\|U_0\|_{\dot{\Gg}^{\s,s}_{\rho_1(0)}}\leq M$, define
\[
t_*:=\sup\left\{t\in[0,T_0]:
\|U(r)\|_{\dot{\Gg}^{\s,s}_{\rho_1(r)}}\leq M
\text{ for every }0\leq r\leq t\right\}.
\]
On $[0,t_*]$, the coefficient in the energy inequality is at least $\eta$ by \eqref{prop:EMHD_main_prop_2_1}, and hence
\begin{equation}\label{eq:base_energy_integrated}
\|U(t)\|_{\dot{\Gg}^{\s,s}_{\rho_1(t)}}^2
+\eta\int_0^t
\|U(r)\|_{\dot{\Gg}^{\s+1,s}_{\rho_1(r)}}^2\,dr
\leq
\|U_0\|_{\dot{\Gg}^{\s,s}_{\rho_1(0)}}^2
\leq M^2.
\end{equation}
Continuity rules out a first crossing above $M$, so $t_*=T_0$. This proves
\eqref{eq:global_induction_statements} for $n=1$.

Assume that \eqref{eq:global_induction_statements} holds for some $n\geq1$.
At time $nT_0$, the restart datum satisfies
\begin{equation}\label{eq:restart_norm_bound}
\|U(nT_0)\|_{\dot{\Gg}^{\s,s}_{\rho_n(nT_0)}}\leq M.
\end{equation}
The time-translated Proposition~\ref{prop:EMHD_main_prop_1}, together with
\eqref{eq:global_barrier_each_radius}, extends the solution from $nT_0$ to
$(n+1)T_0$ at the radius $\rho_n(t)$. The local time is exactly the same $T_0$ because it depends only on the fixed parameters and the common bound $M$. Combining this extension with the induction hypothesis gives
\begin{equation}\label{eq:induction_local_extension}
U\in C\!\left([0,(n+1)T_0];H\cap
\dot{\Gg}^{\s,s}_{\rho_n(t)}\right).
\end{equation}
The positive gap $\k_n$ then yields
\begin{equation}\label{eq:induction_high_regularity}
U\in C\!\left([0,(n+1)T_0];H\cap
\dot{\Gg}^{\s+1,s}_{\rho_{n+1}(t)}\right).
\end{equation}
Applying the energy inequality at $\rho_{n+1}(t)$ and repeating the first-crossing argument leading to
\eqref{eq:base_energy_integrated}, we obtain
\begin{align}
\|U(t)\|_{\dot{\Gg}^{\s,s}_{\rho_{n+1}(t)}}
&\leq
\|U_0\|_{\dot{\Gg}^{\s,s}_{\rho_{n+1}(0)}}
\leq M,
\qquad 0\leq t\leq(n+1)T_0,
\label{eq:induction_bound_next}
\end{align}
and the corresponding norm is nonincreasing. This completes the induction.

We next pass from the radii $\rho_n$ to the limiting smaller radius
$\rho'(t):=\phi(t)+\d'$. Given any finite $T>0$, choose $n$ so large that
$T\leq nT_0$. Since $\rho_n(t)>\rho'(t)$, continuity at the radius
$\rho_n(t)$ implies continuity at $\rho'(t)$ through the bounded multiplier
$e^{-(\d_n-\d')\Lambda^s}$. The induction statements therefore imply
\begin{equation}\label{eq:global_delta_prime}
U\in C\!\left([0,\infty);H\cap
\dot{\Gg}^{\s,s}_{\rho'(t)}\right),
\qquad
\sup_{t\geq0}\|U(t)\|_{\dot{\Gg}^{\s,s}_{\rho'(t)}}\leq M.
\end{equation}
For $0\leq t_1\leq t_2$, choose $n$ so that $t_2\leq nT_0$. Then
\[
\|U(t_2)\|_{\dot{\Gg}^{\s,s}_{\rho_n(t_2)}}
\leq
\|U(t_1)\|_{\dot{\Gg}^{\s,s}_{\rho_n(t_1)}}.
\]
Letting $n\to\infty$ and using dominated convergence in the Fourier sums gives
\begin{equation}\label{eq:delta_prime_monotonicity}
\|U(t_2)\|_{\dot{\Gg}^{\s,s}_{\rho'(t_2)}}
\leq
\|U(t_1)\|_{\dot{\Gg}^{\s,s}_{\rho'(t_1)}}.
\end{equation}

Since \eqref{eq:global_delta_prime} holds for every $\d'<\d$, monotone convergence as
$\d'\uparrow\d$ yields, for each fixed $t\geq0$,
\begin{equation}\label{eq:monotone_convergence}
\|U(t)\|_{\dot{\Gg}^{\s,s}_{\phi(t)+\d}}
=\lim_{\d'\uparrow\d}
\|U(t)\|_{\dot{\Gg}^{\s,s}_{\phi(t)+\d'}}
\leq M.
\end{equation}
This proves membership and the bound at the full radius, but continuity does not follow from monotone convergence alone. Fix an integer $m\geq0$. The datum
$U(mT_0)$ belongs to
$\dot{\Gg}^{\s,s}_{\phi(mT_0)+\d}$ by
\eqref{eq:monotone_convergence} and has norm at most $M$. Applying the time-translated local theorem at the full radius produces a solution in
$C([mT_0,(m+1)T_0];H\cap\dot{\Gg}^{\s,s}_{\phi(t)+\d})$. This solution agrees with the already constructed global solution by uniqueness in any fixed smaller radius
$\phi(t)+\d''$, $\d''<\d$. These intervals cover $[0,\infty)$ and establish full-radius continuity. Finally, letting $\d'\uparrow\d$ in
\eqref{eq:delta_prime_monotonicity} gives monotonicity at the full radius. This proves the proposition.
\end{proof}

\subsection{Proof of Theorem \ref{Thm:main_theorem_2}}

Let $\omega\in\Om^{\a_0,\b}_\mu$. Since $U_0=B_0$, Proposition~\ref{prop:EMHD_main_prop_3} gives a unique global transformed solution in
$C([0,\infty);H\cap\dot{\Gg}^{\s,s}_{\phi(t)+\d})$. Define
$B=e^{\mu W_t\Lambda^s}U$. On the barrier event,
\begin{equation}\label{eq:global_B_bound}
\|B(t)\|_{\dot{\Gg}^{\s,s}_{\d}}
\leq\|U(t)\|_{\dot{\Gg}^{\s,s}_{\phi(t)+\d}}
\leq\|B_0\|_{\dot{\Gg}^{\s,s}_{\a_0+\d}}.
\end{equation}
As in the local proof, strong continuity of
$e^{-(\phi(t)-\mu W_t)\Lambda^s}$ transfers time continuity to $B$, and Proposition~\ref{prop:transform_equivalence} gives the adapted It\^o formulation and uniqueness. Lemma~\ref{lemma:drift_Brownian_motion} yields
\[
\P(\Om^{\a_0,\b}_\mu)
=1-\exp\!\left(-\frac{2\a_0\b}{\mu^2}\right),
\]
which completes the proof.

\section{Conclusion and discussion}\label{sec:conclusion}

In this paper, we studied a stochastic active-vector relaxation of the electron magnetohydrodynamics system without resistivity. The generalized current $\Jalpha B=\curl\Lambda^{\alpha-1}B$ is the periodic power-law specialization of the family introduced by Chae and Jeong \cite{CJ23} and recovers classical EMHD at $\alpha=1$. Multiplicative pseudo-differential noise produces an effective damping term after exponential conjugation. For $3/4<\alpha<1$ and the parameter range stated in Theorem~\ref{Thm:main_theorem_1}, we proved adapted local well-posedness up to a Brownian stopping time and global well-posedness on a quantified high-probability barrier event when the initial Gevrey norm is small relative to the effective noise-induced dissipation. The present estimates do not reach the endpoint $\alpha=1$; treating classical non-resistive EMHD would require exploiting additional structure beyond the convolution bounds used here.

\section*{Author contributions and funding declaration} All authors wrote and reviewed the main manuscript. Qirui Peng is the corresponding author. This work was partially supported by the ONR grant under \#N00014-24-1-2432, the Simons Foundation (MP-TSM-00002783), and the NSF
grant DMS-2420988.

\bibliographystyle{plain}
\bibliography{well_posedness_EMHD}

\end{document}